\documentclass[12pt]{article}
\usepackage{latexsym}
\usepackage{amsfonts} 
 
\newcommand{\proof}{{\noindent \bf Proof. }}
\newtheorem{thm}{Theorem}

\newtheorem{lem}{Lemma}

\begin{document}
\begin{titlepage}
\title{\bf DIVERGENT PERMUTATIONS}
{\author{{\bf Emanuela Fachini}
\\{\tt fachini@di.uniroma1.it}
\\''La Sapienza'' University of Rome
\\ ITALY
\and{\bf J\'anos K\"orner}
\thanks{IASI CNR}
\\{\tt korner@di.uniroma1.it}
\\Rome
\\ ITALY}}

\maketitle
\begin{abstract}
Two permutations of the natural numbers diverge if the absolute value of the difference of their elements in the same position goes to infinity. We show that there exists an infinite number of pairwise divergent permutations of the naturals. We relate this result to more general questions about the permutation capacity of infinite graphs.

\end{abstract}
\end{titlepage}

\section{Introduction} 

We call an arbitrary linear order of the natural numbers an {\em infinite permutation.} Let $G$ be an infinite graph with all the natural numbers as its vertex set. We say that two infinite permutations are {\em $G$-different} if for some integer $i$ they have two adjacent vertices of $G$ in the $i$'th position. Similar concepts have been introduced for permutations of the first $n$ natural numbers in \cite{KM}. Once again, two permutations of $[n]$ are called $G$-different if they put two adjacent vertices of $G$ somewhere in the same position. If the graph $G$ has a finite chromatic number, the largest cardinality of a set of pairwise $G$-different permutations of $[n]$, denoted by  $\omega (G_n)$ grows only exponentially, and its $n$-th root has a limit called the {\em permutation capacity} of the graph $G.$ This quantity is hard to determine even in case of very simple graphs. Let $L$ be the infinite line graph on the set of natural numbers $\mathbb{N}$ where the edges connect consecutive integers $\{i, i+1\}.$ The permutation capacity of this graph,
$$\lim_{n \rightarrow \infty} {1 \over n} \log  \omega(L_n)$$
is known to lie between 0.867 and 1, cf. \cite{KNN} and \cite{sev}. It is conjectured in \cite{KM} that $\omega (L_n)$ equals the middle binomial 
coefficient ${n \choose  \lfloor n/2 \rfloor},$ implying that the upper bound is tight.

It is very natural to consider capacity problems for a class of graphs. The capacity of a graph is an equivalent formulation in  combinatorics of Shannon's 
problem of the capacity of a discrete memoryless channel for probability of error equal to zero \cite{Sh}. One of the breakthrough results of Lov\'asz' 
was the determination of this capacity for the odd cycle of length 5 and all the other self-complementary graphs having a vertex-transitive vertex set \cite{L}. In information theory it is important to extend the capacity problem to a class of channels. 
This means that we are interested in the capacity of an unknown channel, of which we only know to belong to a given class. In information theory this model is called the 
{\em compound channel.} The interested reader can consult \cite{CK} for details. The zero-error version of the compound channel problem was solved by Gargano, K\"orner and Vaccaro \cite{GKV} introducing a novel {\em intertwining} construction technique. As it is explained in \cite{sev}, the intertwining proof technique can be applied to show that the permutation capacity of a finite family of graphs having a finite chromatic number each, is equal to the minimum of the permutation capacities of the graphs. (A channel with minimum permutation capacity is a bottleneck.) The technique 
cannot be applied if the family of graphs involved is infinite. It was shown in \cite{sev} that in some cases of infinite families of digraphs the bottleneck theorem fails to hold.

In this paper we continue the exploration of permutation capacity for infinite classes of graphs in the undirected case.

\section{Divergent permutations}

We partition the edge set of the complete graph with vertex set $\mathbb{N}$ by the graphs $L(k)$ where $L(1)=L$ and the vertices $i$ and $j$ are adjacent in $L(k)$ if the 
absolute value of the difference of the numbers $i$ and $j$ is $k.$ It is proved by K\"orner, Simonyi and Sinaimeri as part of the proof of Proposition 3 in \cite{KSS} that all these graphs have the same permutation capacity. Let us denote by $\cal{L}$ the family of all our graphs. We say that two permutations of $\mathbb{N}$ are  $\cal{L}$-different if they are $L$-different for every $L \in \cal{L}$. We say that two permutations of $\mathbb{N}$
{\em diverge} if the absolute value of the difference between integers in the same position goes to infinity. A seemingly weaker condition for pairs of permutations is to be 
$L$-different for an infinite number of the graphs $L \in \cal{L}$. 

\begin{lem}\label{lem:inf}
There is an infinite number of permutations of $\mathbb{N}$ so that any two of them differ in an infinite number of the graphs $L \in \cal{L}.$
\end{lem}

\proof

We will explicitly construct an infinite family of permutations of $\mathbb{N}$ such that any pair of them will be $L$-different for infinitely many member graphs $L \in \cal{L}.$
For every natural number $i$ we construct an infinite permutation ${\bf x}^i$ as follows. In the $2j$'th position of ${\bf x}^i$ we put $ij.$ This defines the value of ${\bf x}^i$ in even positions. We set the missing natural numbers in the odd positions in the increasing order. Consider now an arbitrary pair, ${\bf x}^i$ and ${\bf x}^k$ of our permutations. 
Suppose that $k>i.$ The difference in the $j$'th even position is $2(k-i)j$ and all these differences are different. We see that $({\bf x}^i_{2j}, {\bf x}^k_{2j})$ form an edge in 
the graph $L(2(k-i)j)$ for every $j \in \mathbb{N}.$

\hfill$\Box$

Further scrutiny of our construction shows that it also proves 

\begin{thm}\label{thm:veg}

There are infinitely many pairwise divergent infinite permutations.  
\end{thm}

\proof

Let us take a closer look at an arbitrary pair of permutations ${\bf x}^i$ and ${\bf x}^k$ from our construction in the proof of the Lemma. If $k>i$, then, obviously, the sequence of the differences $2(k-i)j, j=1, 2, \dots$ goes to infinity. We claim that the sequence of the absolute values of the differences in the odd positions also goes to infinity. This implies that alternating the elements of these two sequences we get a new sequence of differences which also tends to infinity. To verify our statement about the subsequences of the odd positions, 
notice that in the odd positions of ${\bf x}^i $ we find the natural numbers in increasing order, except that in every position which is a multiple of $i$ there is a shift to the left and we see an increase 1 of the value in the position. This means that
$${\bf x}^i_j =j+\Big\lfloor {j \over i}\Big\rfloor   \qquad {\rm and} \qquad      {\bf x}^k_j =j+\Big\lfloor {j \over k }\Big\rfloor        $$
Since $k>i,$ we have 
$${\bf x}^i_j \geq {\bf x}^k_j $$
and
$${\bf x}^i_j - {\bf x}^k_j =\Big\lfloor {j \over i}\Big\rfloor - \Big\lfloor {j \over k }\Big\rfloor .$$
The sequence of these differences is non-negative and goes to infinity. In fact, it is easy to see that the sequence is monotonically increasing, and if 
$j=ikt$ with $t$ tending to infinity, we have
$$\Big\lfloor {j \over i}\Big\rfloor - \Big\lfloor {j \over k }\Big\rfloor=(k-i)t$$
going to infinity.
Hence the sequence of all the differences alternates between two sequences each of which goes to 
infinity, implying the same for the whole sequence.

\hfill$\Box$

This leaves open the question about the cardinality of the largest set of pairwise divergent infinite permutations. 
\section{Related questions}

Infinite permutations can differ in many other interesting ways.
In the already cited paper of K\"orner and Malvenuto \cite{KM} two permutations of $[n]$ were called {\em colliding} if 
in some coordinate they feature two consecutive integers. Denoting by $L$ the infinite line graph over the natural numbers with edges formed by all the pairs of consecutive integers, we now ask whether there is an infinite set of infinite permutations with the property that any pair of them differ in an edge of $L$ in infinitely many coordinates. We will call a pair of infinite permutations with the latter property {\em infinitely colliding}. 

We prove
\begin{thm}\label{thm:cat}
There are infinitely many infinite permutations any pair of which are infinitely colliding.
\end{thm}

\proof

For brevity, we will call an odd integer {\em pure} if it is the $k$'th power of an odd prime with $k>1.$ We will call an integer the {\em successor} of the pure number 
$q$ if it equals $q+1.$ By the Catalan-Mih\u{a}ilescu theorem \cite{M} all these numbers are different. Let $\mathbb{B}$  be the set of all the infinite sequences of zeroes and ones with the exception of the all-zero sequence. We will map the elements of $\mathbb{B}$ injectively to infinite permutations in such a way that any pair of the permutations 
so obtained are infinitely colliding. Let $\bf z \in \mathbb{B}$ be arbitrary. We define the corresponding infinite permutation $\rho$ by an infinite sequence of transpositions in the identity permutation. Notice first that there is at least one $j\in \mathbb{N}$ for which the $j$'th coordinate of $\bf z$  is equal to 1. Whenever this happens, we swap the 
integer $p^j$ in the corresponding position in the identity permutation with its successor, for every pure number, i.e. for every odd prime $p.$ If $\bf z'$ and $\bf z''$ are two different elements of $\mathbb{B}$, then the two permutations differ in an infinite number of swaps so defined. 

\hfill$\Box$

The construction shows that the set of permutations in the proof is uncountable; in fact, it has the cardinality of the set of infinite binary sequences. We observe furthermore, that replacing the graph $L$ by any other member graph of 
$\cal L$, an analogous statement can be proved, since all these graphs contain a semi-infinite path. It is easy to see that the construction works for any infinite graph containing 
an infinite matching, i.e., an infinite set of vertex-disjoint edges. 

To conclude, for an arbitrary infinite graph $G$ with vertex set $\mathbb{N}$ we will say that two infinite permutations ate {\em completely G-different} if in every position their 
corresponding entries form an edge of $G.$ We have

\begin{thm}\label{thm:infin}
There is an infinite number of pairwise completely $K$-different infinite permutations of $K,$ the complete graph with vertex set $\mathbb{N}.$
\end{thm}

\proof

We give an explicit construction. For every natural number $i$ let us divide the sequence of natural numbers into successive disjoint intervals of length $2^i.$ For every such interval let us exchange the second half with the first half. For every $i$ we obtain a different infinite permutation. We claim that any two of these permutations differ in each of their positions. To see this, let us consider the $i$'th and the $j$'th of the permutations with $j>i.$  We look into what happens in the first interval of length $2^j$ of the $j$'th permutation. The first half of the interval in the identity permutation takes the place of the second half. However, in the $i$'th permutation any interval of length $2^i$ in the defining partition keeps its original place as a whole, and changes occur only inside an interval. Hence, in particular, the interval of the first $2^{j-1}$ positions remains, as a whole, in its original place. 
This then means that none of these positions in the $i$'th permutation keeps the same value in both of our permutations. The same reasoning applies to the second half 
of the initial  interval of length $2^j$ in the two permutations. The same relationship holds for the comparison of subsequent intervals. 

\hfill$\Box$

Our result leaves open the question about the largest cardinality of a set of pairwise completely $K$-different permutations. 

It is easy to see that the last theorem can be 
generalized to graphs that are not complete. As an example, consider the infinite graph in which two vertices are adjacent if they have the same parity. For every $i$ apply the previous construction to the odd coordinates and the even coordinates separately. Then, if we consider the $i$'th and the $j$'th permutation, the two are adjacent in every odd 
coordinate for the very same reasoning as before. Hence, the same holds true also for the even coordinates. We immediately realise that our theorem holds, for any natural number $q$ and the infinite graph in which two distinct vertices are adjacent if they are equal modulo $q$. This raises the question of how sparse a graph can be for the last Theorem to remain true.

In conclusion, we briefly recall an intriguing problem about infinite permutations relative to finite graphs from \cite{KMS}. Let $G$ be a finite graph whose vertex set is 
in $\mathbb{N.}$ As before, we will say that two infinite permutations are $G$-different if they have two adjacent vertices of $G$ somewhere in the same position. What is 
the maximum cardinality of a set  of pairwise $G$-different permutations? Clearly, this number, denoted $\kappa(G)$ in \cite{KMS} and \cite{sev}, is finite and can be considered an interesting parameter of graph $G.$ Its value is unknown even for complete graphs on at least 4 vertices.

\end{document}